\documentstyle[12pt]{article}       

\begin{document}

\newcommand{\comment}[1]{}    
\newcommand{\hs}{\enspace}
\newcommand{\hhs}{\thinspace}
\newcommand{\real}{\ifmmode {\rm R} \else ${\rm R}$ \fi}
\def\nat{\hbox{\vrule height 7pt width .7pt depth 0pt\hskip -.5pt\bf N}}
\newcommand{\qed}{\hfill{\setlength{\fboxsep}{0pt}
                  \framebox[7pt]{\rule{0pt}{7pt}}} \newline}
\newcommand{\eqed}{\qquad{\setlength{\fboxsep}{0pt}
                  \framebox[7pt]{\rule{0pt}{7pt}}} }
\newtheorem{theorem}{Theorem}
\newtheorem{lemma}[theorem]{Lemma}         
\newtheorem{corollary}[theorem]{Corollary}
\newtheorem{definition}[theorem]{Definition}
\newtheorem{claim}[theorem]{Claim}
\newtheorem{conjecture}[theorem]{Conjecture}


\newcommand{\proof }{{\bf Proof: \enspace}}          


\def\edge{\leftrightarrow}
\def\noedge{\not\leftrightarrow}
\def\twoedge{\Leftrightarrow}


\title{Realizing Degree Imbalances in Directed Graphs}
\author{Dhruv Mubayi \ \ and \ Todd G. Will\\
        Davidson College, Davidson NC 28036 \\
	\\
	 Douglas B. West \\
	University of Illinois, Urbana, IL 61801}
\date{}
\maketitle

\begin{abstract}
In a directed graph, the imbalance of a vertex is its outdegree minus its
indegree.  We characterize the sequences that are realizable as the sequence of
imbalances of a simple directed graph.  Moreover, a realization of a realizable
sequence can be produced by a greedy algorithm. 
\end{abstract}


A sequence of integers is {\em graphic} if there is a simple undirected graph
for which it is the degree sequence.  There are well-known characterizations
of graphic sequences, both recursive and non-recursive.  We obtain analogous
results for the ``imbalance sequences'' of digraphs with no repeated arcs.  We
define the {\em imbalance} of a vertex in a digraph to be the difference between
the number of exiting arcs and the number of entering arcs.  Symbolically,
$b(v)=d^{+}(v)-d^{-}(v)$, where $d^{+}(v)$ is the out-degree of $v$ and
$d^{-}(v)$ is the in-degree.  A sequence $b = b_1,b_2,\ldots,b_n,$ is {\em
realizable} if there exists a simple digraph $G$ with vertices
$v_1,v_2,\ldots ,v_n$ such that $b(v_i)=b_i$.

If we allow repeated arcs, then the trivial necessary condition $\sum b_i=0$
is sufficient, using only arcs from vertices with positive imbalance to
vertices with negative imbalance; this is analogous to the observation that
when $\sum d_i$ is even, the sequence $d$ of nonnegative integers is the degree
sequence of some undirected graph, allowing loops and multiple edges.  Hence we
forbid repeated arcs.  It does not matter whether we forbid loops or pairs of
opposed arcs, since these have no affect on the imbalance sequence.

Our results are analogous to the known conditions for graphic sequences.
Havel \cite{Hav} and Hakimi \cite{Hak} independently showed that a
non-increasing sequence of nonnegative integers $d = d_1, d_2,\ldots,d_n$ is
graphic if and only if $\hat d= \hat d_2,\ldots, \hat d_n$ is graphic, where
$\hat d$ is formed from $d$ by deleting $d_1$ and subtracting $1$ from the
$d_1$ largest remaining elements of $d$.  In this context, we define $d_i$ to be
{\em larger} than $d_j$ if $d_i > d_j$ or if $d_i = d_j$ and $i > j$, so that
$\hat d$ is also non-increasing.  This result provides a recursive procedure
that tests whether a sequence is graphic.  If so, retracing the computation
produces a realization; if $\hat G$ realizes $\hat d$, then adding a vertex
whose neighborhood corresponds to the $d_1$ largest entries of $\hat d$ forms a
graph $G$ that realizes $d$.

Erd\H{o}s and Gallai \cite{EG} provided a non-recursive characterization of
graphic sequences: the nonnegative integer sequence $d = d_1, d_2,\ldots,d_n$ 
is graphic if and only if $\sum d_i$ is even and
$\sum_{i=1}^k d_i\leq k(k-1)+\sum_{j=k+1}^n \min \{k, d_j\} $
for all $1\leq k\leq n$.  The inequalities are necessary because from the
set $S$ of $k$ vertices with highest degree there are at most
$\sum_{j=k+1}^n\min\{k,d_j\}$ edges to $V(G)-S$, and there are
at most $k(k-1)/2$ edges within $S$.  Recently, Aigner and Triesch \cite{AT}
presented a short proof of sufficiency using ideals in the dominance order on
sequences with fixed sum.

A greedy algorithm for realizing an imbalance sequence directs arcs from the
vertex with greatest imbalance ($b_1$) to the $b_1$ vertices with smallest
imbalance.  This heuristic is natural in the sense that the vertices with large
(positive) imbalances require high out-degree, and vertices with small (i.e.\
negative) imbalances require high in-degree.  This suggests a modified sequence
$\hat b$ and a recursive test, analogous to the Havel-Hakimi test.

We also prove a non-recursive characterization analogous to the Erd\H{o}s-Gallai
result.  If $b=b_1,b_2,\ldots,b_n$ is realizable, then $\sum b_i=0$, because
each edge contributes positively and negatively to the total imbalance.
Also, $\sum_{i=1}^k b_i\leq k(n-k)$; letting $S$ be the set of $k$ vertices
with largest imbalance, the edges within $S$ contribute nothing to
$\sum_{i=1}^k b_i$, and the pairs $S\times (V(G)-S)$ contribute at most one
each.

These obvious necessary conditions are also sufficient; we present two proofs.
The first shows that the greedy algorithm suggested above produces a realization
when one exists.  The second proof, based on the method of ideals used by Aigner
and Triesch, is much shorter but does not yield such a fast construction
algorithm.

\begin{definition}
A non-increasing sequence of integers $a=a_1,a_2,\ldots ,a_n$ is a {\em feasible
sequence} if the sum of the elements is zero and $\sum_{i=1}^ka_i\leq k(n-k)$
for all $1\leq k\leq n$.
\end{definition}

From a feasible sequence $a=a_1,a_2,\ldots ,a_n$, we form an associated
sequence $\hat a=\hat a_2,\hat a_3,\ldots ,\hat a_n$ by deleting $a_1$ and
adding $1$ to the $a_1$ smallest elements of $a$, where $a_i$ is smaller than 
$a_j$ if $a_i<a_j$ or if $a_i=a_j$ and $i<j$.  The sequence $\hat a$ is the
result of applying one step of the greedy algorithm to the sequence $a$.
The complication in determining which entries of $a$ are augmented to obtain
$\hat a$ occurs when $a_{n-a_1+1}=a_{n-a_1}$; in this case there is a gap 
consisting of elements to which we add zero. We endure this complication to
ensure that $\hat a$ is also non-increasing.  The following example
with $n=9$ and $a_1=5$ produces such a gap, since $a_5=a_4=2$.
$$
\begin{array}{c}
a: \\  
\\ 
\hat a: 
\end{array}
\begin{array}{ccccccccc}
5 & 3 & 2 & 2 & 2 & 2 & -5 & -5 & -6 \\ 
\cdot & 0 & 1 & 1 & 0 & 0 & 1 & 1 & 1 \\ 
\cdot & 3 & 3 & 3 & 2 & 2 & -4 & -4 & -5 
\end{array}
$$
The values of $a_i$ in the gap and the values of $a_i$ to the left of the gap
are all equal.  When this occurs in a feasible sequence, a stronger statement
can be made about the partial sums.

\begin{lemma}
If $a$ is a feasible sequence and $a_k=a_{k+1}=\cdots=a_{k+m}$, then
   $\sum_{i=1}^k a_i \le k(n-k)-m$.
\end{lemma}
\proof
If $a_k\le n-2k-m$, then
$$\sum_{i=1}^k a_i = (\sum_{i=1}^{k-1}  a_i)+a_k \le (k-1)(n-k+1)+(n-2k-m) = k(n-k)-m-1.$$
If $a_k > n-2k-m$, then
$$\sum_{i=1}^{k} a_i = (\sum_{i=1}^{k+m}  a_i)-ma_k\le (k+m)(n-k-m)-m(n-2k-m+1) =k(n-k)-m.\eqed$$

\begin{theorem}  If $a$ is feasible, then $\hat a$ is feasible. \end{theorem}
\proof
Because $\sum_{i=1}^{n}  a_i=0$ and $a_1\le 1\cdot(n-1)=n-1$, the sequence
$\hat a$ is well defined.  By the definition of the $a_1$ ``smallest'' elements,
$\hat a$ is non-increasing.  The construction of $\hat a$ distributes $a_1$
among the other entries, so the sum is still 0.  It remains to verify the
condition on partial sums, which is $\sum_{i=2}^{k}  \hat a_i\le(k-1)(n-k)$ for 
$2\le k\le n$.

Suppose $a_1=n-r$, and suppose the construction of $\hat a$ involves a gap
of $s$ entries to which we add 0, where $s\ge 0$.  Beginning with the leftmost
entry to which 1 is added, there are $(n-r)+s$ positions in the sequence;
hence $s<r$ and $a_{r-s+1}$ is the leftmost entry to which $1$ is added.
If $k\le r-s$, then
$$\sum_{i=2}^k  \hat a_i = \sum_{i=2}^k a_i \le (k-1)(n-r) \le (k-1)(n-k).$$

Now suppose $k>r-s$.  Let $\alpha=\sum_{i=2}^{k} \hat a_i -\sum_{i=2}^{k} a_i$;
this is the number of added 1's in the sum.  Let $t$ be the number of positions
in the sequence that are to the left of the gap and receive an augmentation of
$1$ ($t=0$ if $s=0$).  The rightmost entry to which $0$ is added is $a_{r+t}$.
$$
\begin{array}{ccccc}
a: 	& a_1 	& a_2 		& \cdots 	& a_{r-s}\\  
	&   	& 0 		& \cdots 	& 0 \\ 
\hat a: 	& 	& \hat a_2 	& \cdots	& \hat a_{r-s}
\end{array}
\stackrel{t}{\overbrace{
\begin{array}{ccc}
a_{r-s+1}   	& \cdots	& a_{r-s+t} 		\\
1		&  \cdots	& 1  		  	\\
\hat a_{r-s+1} &  \cdots	& \hat a_{r-s+t} 	\\
\end{array}
}}
\stackrel{s}{\overbrace{
\begin{array}{ccc}
a_{r-s+t+1}		& \cdots	& a_{r+t}		\\
0			& \cdots	& 0			\\
\hat a_{r-s+t+1}	& \cdots	& \hat a_{r+t}		\\
\end{array}
}} 
\begin{array}{ccc}
a_{r+t+1}		& \cdots	& a_{n}		\\
1			& \cdots	& 1			\\
\hat a_{r+t+1}		& \cdots	& \hat a_{n}		
\end{array}
$$

If $r-s<k< r+t$, then $a_k=\cdots=a_{r+t}$.  By the Lemma,
$\sum_{i=1}^{k}  a_i \le k(n-k)-(r+t-k)$.  In this case, $\alpha\le t$.
If $k\ge r+t$, then $\sum_{i=1}^{k}  a_i\le k(n-k)$ and $\alpha =k-r$.
In either case, we have $(\sum_{i=1}^{k}  a_i) +\alpha\le k(n-k)+(k-r)$.  Hence
$$\sum_{i=1}^{k} \hat a_i=(\sum_{i=1}^{k}  a_i)+\alpha-a_1\le k(n-k)+(k-r)-(n-r)=(k-1)(n-k).\eqed$$

A simple inductive argument now completes the proof of the sufficiency of the
condition and shows that the greedy algorithm produces a realization if one
exists.  We also present a second proof, using the sufficiency method of Aigner
and Triesch [4].  The idea here is to define a partial order $P$ on sequences
such that 1) the realizable sequences form an ideal (down-set), and 2) the
maximal sequences in $P$ that satisfy the desired condition (here feasibility)
are realizable.  This implies that all sequences satisfying the condition are
realizable, because they are dominated in $P$ by a realizable sequence.

\begin{corollary}
A sequence $a$ is realizable as an imbalance sequence if and only if $a$ is
feasible.
\end{corollary}

\proof 
{\bf 1:}
We have already argued the necessity of the condition.  For sufficiency,
we used induction on the length $n$ of $a$.  If $n=1$, then $a_1=0$ and $K_1$
has imbalance sequence $a$.  If $n>1$, then Theorem 3 implies that $\hat a$ is a
feasible sequence of length $n-1$.  By the induction hypothesis, there is a
graph $G'$ with vertices $v_2,v_3,\ldots ,v_n$ such that $b(v_i)=\hat a_i$.
Form $G$ by adding the vertex $v_1$ and arcs from $v_1$ to $v_j$ for each
$j\in V(G')$ such that $\hat a_j=a_j+1$.  The graph $G$ realizes $a$.

\proof
{\bf 2:}
Consider the domination order $P$ on non-increasing sequences of integers
with sum $0$; the order is defined by $a\ge b$ if
$(\sum_{i=1}^{k}  a_i)\ge(\sum_{i=1}^{k}  b_i)$ for all $k$.
In the subposet of $P$ induced by the feasible sequences, we claim that only
the sequence $n-1,n-3,\ldots,-(n-3),-(n-1)$ is maximal.  It achieves each
constraint with equality, and therefore by the definition of the order it
dominates all other feasible sequences.  Furthermore, this is the imbalance
sequence of the transitive tournament with $n$ vertices and hence is realizable.

It remains to show that the realizable sequences in $P$ form an ideal.
If $b>b'$ in $P$, then $b'$ is obtained from $b$ by a sequence of unit shifts in
which some position $i$ decreases by $1$ and some later position $j$ increases
by one; it suffices to show that every such shift maintains realizability.
Suppose $b$ is realizable, and consider positions $i,j$ for such a shift.
Because $b_i>b_j$, there is some other vertex $z$ in the realization such that 
     $i \to z \noedge j$    or    $i \to z \to j$    or   $ i \noedge z \to j$.
In these three cases, respectively, we transform the graph to
    $i \noedge z \leftarrow j$    or   $i \noedge z \noedge j$   or    $i \leftarrow z \noedge j$.
The only change in the sequence is that $b_i$ goes down by $1$ and $b_j$ goes up
by $1$, so we have realized $b'$ also as an imbalance sequence.  \qed

\end{document}